\documentclass[article]{elsarticle}

\usepackage{hyperref}
\usepackage{geometry}
\usepackage{amsmath,amssymb,enumerate,amsfonts}

\usepackage[normalem]{ulem} 
\usepackage{longtable,threeparttable,float}
\usepackage[table]{xcolor}
\usepackage{graphicx, subfigure}
\usepackage{caption}
\usepackage{algorithm}%
\usepackage{algorithmicx}%
\usepackage{algpseudocode}%
\usepackage{epstopdf}
\usepackage{tabularx}
\usepackage[justification=centering]{caption}
\usepackage {color}
\usepackage{hyperref} 
\hypersetup{
	colorlinks=true,
	citecolor=blue,
	linkcolor=blue,
	urlcolor=green}

\usepackage{setspace}
\newenvironment{newlist}[1]%
{\begin{list}{}{\settowidth{\labelwidth}{\bf #1}%
			\setlength{\leftmargin}{\labelwidth}%
			\addtolength{\leftmargin}{\labelsep}%
			}}%
	{\end{list}}

\usepackage[figuresright]{rotating}

\journal{Journal of \LaTeX\ Templates}

\newtheorem{definition}{Definition}[section]
\newtheorem{remark}{Remark}[section]
\newtheorem{example}{Example}[section]
\newtheorem{theorem}{Theorem}[section]
\newtheorem{proposition}{Proposition}[section]
\newtheorem{lemma}{Lemma}[section]

\begin{document}

\begin{frontmatter}

\title{On the convergence of conditional gradient method for unbounded multiobjective optimization problems}

\author[mymainaddress,mysecondaryaddress]{Wang Chen}
\ead{chenwangff@163.com}

\author[mythirddaryaddress]{Yong Zhao}	
\ead{zhaoyongty@126.com}

\author[mymainaddress]{Liping Tang}
\ead{tanglps@163.com}

\author[mymainaddress]{Xinmin Yang}
\ead{xmyang@cqnu.edu.cn}

\address[mymainaddress]{National Center for Applied Mathematics in Chongqing, Chongqing Normal University, Chongqing, 401331, China}
\address[mysecondaryaddress]{School of Mathematical Sciences, University of Electronic Science and Technology of China, Chengdu, 611731, China}
\address[mythirddaryaddress]{College of Mathematics and Statistics, Chongqing Jiaotong University, Chongqing, 400074, China}

\begin{abstract}
This paper focuses on developing a conditional gradient algorithm for multiobjective optimization problems with an unbounded feasible region. We employ the concept of recession cone to establish the well-defined nature of the algorithm. The asymptotic convergence property and the iteration-complexity bound are established under mild assumptions. Numerical examples are provided to verify the algorithmic performance.
\end{abstract}


\end{frontmatter}

\section{Introduction}

Multiobjective optimization refers to the problem of optimizing several objective
functions simultaneously. These problems often entail trade-offs between conflicting and competing objectives. For instance, designing a car may involve concurrently optimizing fuel efficiency, safety, comfort, and aesthetics. This type of problem has applications in  engineering \cite{R_m2013}, finance \cite{Z_m2015}, environmental analysis \cite{F_o2001}, management science \cite{T_a2010}, machine learning \cite{J_m2006,sener2018}, etc. 

The multiobjective optimization problem has the following form:
\begin{equation}\label{mop}
	\text{min}\quad F(x)~~~~~
	\text{s.t.}\quad x\in \Omega,
\end{equation}
where $F(x) = (F_{1}(x),F_{2}(x),...,F_{m}(x))$ is a vector-valued function with each $F_{i}$ being  continuously differentiable, and $\Omega\subset\mathbb{R}^{n}$ is a feasible region. When $\Omega=\mathbb{R}^{n}$, numerous descent algorithms are currently developed to solve \eqref{mop}; see, for example,  \cite{fliege2000steepest,fliege2009newton,lucambio2018nonlinear,lapucci2023limited}. In scenarios where $\Omega$ is assumed to be a compact set  (i.e., bounded and closed) and convex set, the conditional gradient methods \cite{assunccao2021conditional,chen2023conditional} have been devised for solving \eqref{mop}. In many practical applications, however, the feasible region $\Omega$ may be unbounded, which limits the applicability of the conditional gradient methods. Some motivating examples can be found in the multiobjective optimization literature  \cite{lin2005on,hoa2007unbounded,li2008equivalence,wagner2023algorithms,huong2020geoffrion,kov2022convex,meng2022portfolio}. The major contribution of this paper is to generalize the traditional conditional gradient method \cite{assunccao2021conditional,chen2023conditional} to solve \eqref{mop} with computational guarantees, where $\Omega$ is nonempty closed and convex (not necessarily compact).

The rest of the work is organized as follows. Section
\ref{sec:2} provides some basic definitions, notations and auxiliary results. Section \ref{sec:3} gives the conditional gradient algorithm. Section \ref{sec:4} is devoted to the investigation of the convergence properties. Section \ref{sec:5} includes numerical experiments to demonstrate the algorithm's performance.

\section{Preliminaries}\label{sec:2}

Denote by $\langle \cdot,\cdot\rangle$ and $\|\cdot\|$, respectively, the usual inner product and the norm in $\mathbb{R}^{n}$. Let  $\langle m\rangle=\{1,2,\ldots,m\}$ and $e=(1,1,\ldots,1)^{\top}$. Recall that the dual cone  of a cone $C$ in $\mathbb{R}^{n}$ and its interior are, respectively, defined by
$C^{*} = \{ y^{\ast} \in \mathbb{R}^{n} : \langle y,y^{\ast}\rangle \geq 0 ~{\rm for~all} ~y\in C\}$
and
\begin{equation}\label{int}
		{\rm int}(C^{*}) = \{y^{\ast} \in\mathbb{R}^{n}: \langle y,y^{\ast}\rangle > 0,\forall y \in C \setminus\{ 0 \} \}.
\end{equation}
 For any given nonempty set $A\subset \mathbb{R}^{n}$, we define the recession cone of  $A$ (see \cite[pp. 81]{rockafellar1998}), denoted by $A^{\infty}$, as
\begin{equation*}\label{recession_cone1}
	A^{\infty}=\left\{d\in\mathbb{R}^{n}:\exists \{x^{k}\}\subset A~{\rm and}~\exists \{\lambda_{k}\}{~\rm with}~\lambda^{k}\downarrow0~{\rm such~that}~\lim_{k\rightarrow\infty}\lambda_{k}x^{k}=d\right\}.
\end{equation*}
When $A$ is closed and convex, its recession cone can be determined by the following formula:
\begin{equation}\label{recession_cone2}
	A^{\infty}=\{d\in\mathbb{R}^{n}:x+td\in A,\forall x\in A, t\geq0\}.
\end{equation}
The importance of the recession cone is revealed by the key property that $A$ is bounded if and only if $A^{\infty}=\{0\}$ (see \cite[pp. 81]{rockafellar1998}).

Let $\mathbb{R}_{+}^{m}$ and $\mathbb{R}_{++}^{m}$ denote the non-negative orthant and positive orthant of $\mathbb{R}^{n}$, respectively. We may consider the partial order $\preceq~(\prec)$ induced by $\mathbb{R}_{+}^{m}~ (\mathbb{R}_{++}^{m})$: for any $x,y\in\mathbb{R}^{m}$, $x\preceq  y~(x\prec  y)$ if and only if $y-x\in\mathbb{R}_{+}^{m}~ (y-x\in\mathbb{R}_{++}^{m})$. The Jacobian of $F$ at $x=(x_{1},x_{2},\ldots,x_{n})\in\mathbb{R}^{n}$ is denoted by $JF(x)=[\nabla F_{1}(x)~\nabla F_{2}(x)~\ldots~\nabla F_{m}(x)]^{\top}$. Recall that $F$ is convex on $\Omega$ if and only if $JF(y)(x-y)\preceq F(x)-F(y)$ for all $x,y\in\Omega$ and all $\lambda\in[0,1]$ 
 (see \cite{J2011}).

A point $\bar{x}\in\Omega$ is called a Pareto optimal solution of (\ref{mop}) if there does not exist any other $x\in\Omega$ such that $F(x)\preceq F(\bar{x})$ and $F(x)\neq F(\bar{x})$, and a point $\bar{x}\in\Omega$ is called a weak Pareto optimal solution of (\ref{mop}) if there does not exist any other $x\in\Omega$ such that $F(x)\prec F(\bar{x})$ (see \cite{miettinen1999nonlinear}). A necessary, but not sufficient, first-order optimality condition for  (\ref{mop}) at $\bar{x}\in\Omega$, is
\begin{equation}\label{fir_ord_opt}
	JF(\bar{x})(\Omega-\bar{x})\cap(-\mathbb{R}_{++}^{m})=\emptyset,
\end{equation}
where $JF(\bar{x})(\Omega-\bar{x})=\{JF(\bar{x})(u-\bar{x}):u\in\Omega\}$ and
$$JF(\bar{x})(u-\bar{x})=(\langle\nabla F_{1}(\bar{x}),u-x\rangle,\langle\nabla F_{2}(\bar{x}),u-\bar{x}\rangle,\ldots,\langle\nabla F_{m}(\bar{x}),u-\bar{x}\rangle)^{\top}.$$

\begin{definition}\normalfont\label{pare_sta}
	A point $\bar{x}\in\Omega$ satisfying (\ref{fir_ord_opt}) is called a Pareto critical point of  (\ref{mop}).
\end{definition}

\begin{remark}\normalfont
	As mentioned in \cite{assunccao2021conditional}, the geometric optimality condition \eqref{fir_ord_opt} can also be equivalently expressed as
	\begin{equation}\label{opt_con}
		\max_{i\in\langle m\rangle}\{\langle \nabla F_{i}(\bar{x}),u-\bar{x}\rangle\}\geq0\quad {\rm for~ all}~u\in\Omega.
	\end{equation}
\end{remark}

\begin{lemma}\label{global}{\rm\cite{assunccao2021conditional}}
	If $F$ is convex on $\Omega$ and $\bar{x}\in\Omega$ is a Pareto critical point, then $\bar{x}$ is also a weak Pareto 
	optimal solution of \eqref{mop}.
\end{lemma}

\begin{lemma}{\rm\cite{beck2017first}}\label{sub_rate}
	Let $\{a_{k}\}$ be a sequence of nonnegative real numbers satisfying for any $k \geq 0$,
	$a_{k}-a_{k+1}\geq a_{k}^{2}/\gamma $
	for some $\gamma>0$. Then, for any $k\geq1$,
	$a_{k}\leq\gamma/k.$
\end{lemma}

We end this section by assuming each gradient function $\nabla F_{i}$ is Lipschitz continuous with Lipschitz constant $L_{i}>0$ on $\Omega$, i.e., 
$
\|\nabla F_{i}(x)- \nabla F_{i}(y)\|\leq L_{i}\|x-y\|
$
for all $x,y\in\Omega$ and $i\in\langle m\rangle$. In the paper, let $L=\max_{i\in\langle m\rangle}L_{i}$. 

\section{The conditional gradient algorithm}\label{sec:3}

Given $x \in\Omega$, we consider the following auxiliary scalar optimization problem:
\begin{equation}\label{sca_pro}
	\min_{u\in\Omega}\max_{i\in\langle m\rangle}\{\langle \nabla F_{i}(x),u-x\rangle\}.
\end{equation}
Note that the existence of solution for  \eqref{sca_pro} cannot be guaranteed since $\Omega$ is not assumed to be bounded. Listed below is a mild yet key assumption regarding each gradient function, which will be used to show the sequence $\{x^{k}\}$ produced by the conditional gradient algorithm is well-defined.

\begin{description}
	\item[{\rm\bf(A1)}] Each gradient function $\nabla F_{i}$ satisfies
	$
		\nabla F_{i}(x)\in{\rm int}(\Omega^{\infty})^{*}
	$
	for all $x\in\Omega$ and $i\in\langle m\rangle$.
\end{description}

\begin{remark}\normalfont
	Assumption (A1) holds trivially whenever the closed convex set $\Omega$ is bounded. Indeed, $\Omega$ is bounded if and only if $\Omega^{\infty}=\{0\}$, and thus ${\rm int}(\Omega^{\infty})^{*} =\mathbb{R}^{n}$.
\end{remark}

Next, under (A1), we present some results that guarantee the existence of solution of \eqref{sca_pro}.

\begin{proposition}\label{prop1}
	Assume that (A1) holds. For all $x\in\Omega$, the set 
	\begin{equation*}
		\Omega_{1}(x)=\left\{u\in\Omega:\max_{i\in\langle m\rangle}\{\langle\nabla F_{i}(x),u-x\rangle\}\leq 0\right\} 
	\end{equation*}
	is compact. Furthermore, the problem \eqref{sca_pro} has a solution.
\end{proposition}
\textbf{Proof.}
	It follows from (A1) and \eqref{int} that $\langle\nabla F_{i}(x),d\rangle>0$ for any $d\in\Omega^{\infty}\backslash\{0\}$ and $i\in\langle m\rangle$. This implies that
	\begin{equation}\label{pro1}
		\max_{i\in\langle m\rangle}\{\langle\nabla F_{i}(x),d\rangle\}>0
	\end{equation}
	for all $d\in\Omega^{\infty}\backslash\{0\}$. Assume by contradiction that $\Omega_{1}(x)$ is unbounded. Therefore, there exists a sequence $\{u^{k}\}\subset \Omega_{1}(x)$ such that $\lim_{k\rightarrow\infty}\|u^{k}\|=\infty$. Define $\lambda_{k}=1/\|u^{k}\|$. Then, we have $\lim_{k\rightarrow\infty}\lambda^{k}=0$. Clearly, for all $k\geq 0$, 
	$
	\|\lambda_{k}u^{k}\|=\|u^{k}/\|u^{k}\|\|=1.
	$
	This means that there exist subsequences $\{u^{k_{j}}\}\subset \Omega_{1}(x)$ and $\{\lambda_{k
		_{j}}\}\subset(0,\infty)$ with $\lim_{j\rightarrow\infty}\lambda_{k_{j}}=0$ such that 
	\begin{equation}\label{pro2}
		\lim_{j\rightarrow\infty}\lambda_{k_{j}}u^{k_{j}}=\bar{d}\in\Omega^{\infty}.
	\end{equation}
	From the definition of $ \Omega_{1}(x)$ and the positiveness of $\lambda_{k
		_{j}}$, we have
	\begin{equation*}
		0\geq \lambda_{k_{j}}\max_{i\in\langle m\rangle}\{\langle\nabla F_{i}(x),u^{k_{j}}-x\rangle\}\geq\max_{i\in\langle m\rangle}\{\langle\nabla F_{i}(x),\lambda_{k_{j}}u^{k_{j}}\rangle\}-\lambda_{k_{j}}\max_{i\in\langle m\rangle}\{\langle\nabla F_{i}(x),x\rangle\}.
	\end{equation*}
	Taking the limit as $j\rightarrow\infty$ in the above relation, and observing \eqref{pro2}, we obtain 
	$\max_{i\in\langle m\rangle}\{\langle\nabla F_{i}(x),\bar{d}\rangle\}\leq0,$
	contradicting \eqref{pro1} and concluding the proof.
\qed

\begin{proposition}\label{prop2}
	Assume that (A1) holds. If $\Omega_{2}\subset \Omega$ is a bounded set, then the set
	\begin{equation}\label{prop2_0}
		\bigcup_{x\in\Omega_{2}}\left\{p(x)\in\Omega:p(x)\in\mathop{\rm argmin}_{u\in\Omega}\max_{i\in\langle m\rangle}\{\langle\nabla F_{i}(x),u-x\rangle\}\right\} 
	\end{equation}
	is bounded.
\end{proposition}
\textbf{Proof.}
	Assume by contradiction that the set in \eqref{prop2_0} is unbounded. Then, there exists $\{x^{k}\}\subset \Omega_{2}$ and $\{p(x^{k})\}\subset \Omega$ such that $\lim_{k\rightarrow\infty}\|p(x^{k})\|=\infty$. Let $\lambda_{k}=1/\|p(x^{k})-x^{k}\|$. Then, $\lim_{k\rightarrow\infty}\lambda_{k}=0$ because $\Omega_{2}$ is bounded. Clearly, for all $k\geq 0$, we get
	$\|\lambda_{k}(p(x^{k})-x^{k})\|=\|(p(x^{k})-x^{k})/\|p(x^{k})-x^{k}\|\|=1$, which implies that there exist subsequences $\{x^{k_{j}}\}\subset \Omega_{2}$,  $\{p(x^{k_j})\}\subset \Omega$ and $\{\lambda_{k
		_{j}}\}\subset(0,\infty)$ such that
		
		$$\lim_{j\rightarrow\infty}x^{k_{j}}=\bar{x}\quad{\rm and}\quad\lim_{j\rightarrow\infty}\lambda_{k_{j}}(p(x^{k_j})-x^{k_{j}})=\bar{d}.$$
	
	Since $\{x^{k}\}\subset\Omega_{2}\subset\Omega$, $\{p(x^{k})\}\subset\Omega$ and $\Omega$ is a convex set, we have
	$
		x^{k}+\alpha(p(x^{k})-x^{k})\in\Omega
	$
	for all $\alpha\in(0,1)$. Therefore, 
	\begin{equation*}
		\begin{aligned}
			\lim_{j\rightarrow\infty}\lambda_{k_{j}}(x^{k_{j}}+\alpha(p(x^{k_{j}})-x^{k_{j}}))
			&=\lim_{j\rightarrow\infty}(\lambda_{k_{j}}x^{k_{j}}+\alpha\lambda_{k_{j}}(p(x^{k_{j}})-x^{k_{j}}))\\
			&=\lim_{j\rightarrow\infty}\lambda_{k_{j}}x^{k_{j}}+\alpha\lim_{j\rightarrow\infty}\lambda_{k_{j}}(p(x^{k_{j}})-x^{k_{j}})\\
			&=\alpha \bar{d}\in\Omega^{\infty},
		\end{aligned}
	\end{equation*}
	and thus $\bar{d}\in\Omega^{\infty}$ because $\Omega^{\infty}$ is a cone. By (A1), for all $x\in\Omega$ and $i\in\langle m\rangle$, we get 
	\begin{equation}\label{prop2_1}
		\langle\nabla F_{i}(x),\bar{d}\rangle>0.
	\end{equation}
	From \eqref{prop2_0}, we get $p(x^{k_{j}})\in\mathop{\rm argmin}_{u\in\Omega}\max_{i\in\langle m\rangle}\{\langle\nabla F_{i}(x^{k_{j}}),u-x^{k_{j}}\rangle\}$, and observing that $\{x^{k_{j}}\}\subset\Omega_{2}\subset\Omega$, it holds that
	\begin{equation}\label{prop2_2}
		\max_{i\in\langle m\rangle}\{\langle\nabla F_{i}(x^{k_{j}}),p(x^{k_{j}})-x^{k_{j}}\rangle\}\leq\max_{i\in\langle m\rangle}\{\langle\nabla F_{i}(x^{k_{j}}),x^{k_{j}}-x^{k_{j}}\rangle\}=0.
	\end{equation}
	Owing to $\{\lambda_{k_{j}}\}\subset(0,\infty)$, \eqref{prop2_2} implies that
	$
		\max_{i\in\langle m\rangle}\{\langle\nabla F_{i}(x^{k_{j}}),\lambda_{k_{j}}(p(x^{k_{j}})-x^{k_{j}})\rangle\}\leq0,
	$
	i.e., 
	$$
		\langle\nabla F_{i}(x^{k_{j}}),\lambda_{k_{j}}(p(x^{k_{j}})-x^{k_{j}})\rangle\leq0
	$$
	for all $i\in\langle m\rangle$.
	Taking the limit as $j\rightarrow\infty$ in the above relation, we have
	$
		\langle\nabla F_{i}(\bar{x}),\bar{d}\rangle\leq0
	$
	for all $i\in\langle m\rangle$, which is a contradiction to \eqref{prop2_1}. Thus, the proof is complete.
\qed

Denote by $p(x)$ the optimal solution of  \eqref{sca_pro}, i.e.,
\begin{equation}\label{opt_sol}
	p(x)\in\mathop{\rm argmin}_{u\in\Omega}\max_{i\in\langle m\rangle}\{\langle \nabla F_{i}(x),u-x\rangle\}.
\end{equation}
According to Propositions \ref{prop1} and \ref{prop2}, $p(x)$ is well-defined. The optimal value of  \eqref{sca_pro} is denoted by $\theta(x)$, i.e.,
\begin{equation}\label{theta_xk}
	\theta(x)=\max_{i\in\langle m\rangle}\{\langle\nabla F_{i}(x),p(x)-x\rangle\}.
\end{equation}

\begin{lemma}{\rm \cite{assunccao2021conditional}}
	Let $\theta:\Omega\rightarrow\mathbb{R}$ be as in \eqref{theta_xk}. Then, 
	\begin{enumerate}[{\rm (i)}]\setlength{\itemsep}{-0.03in}
		\item $\theta(x)\leq0$ for all $x\in\Omega$;
		\item $\theta(x)=0$ if and only if $x\in\Omega$ is a Pareto critical point.
	\end{enumerate}
\end{lemma}

The general scheme of the conditional gradient (CondG) algorithm for solving  \eqref{mop} is summarized as follows.

\begin{newlist}{Step1: }
	\setlength{\itemsep}{-0.03in}
	\item[CondG algorithm.]
	\item[Step 0] Choose $x^{0}\in\Omega$. Compute $p(x^{0})$ and $\theta(x^{0})$ and initialize $k\leftarrow0$.
	\item[Step 1] If $\theta(x^{k})=0$, then \textbf{stop}.
	\item[Step 2] Compute $d(x^{k})=p(x^{k})-x^{k}$.
	\item[Step 3] Compute the step size $t_{k}\in(0,1]$ by a step size strategy and set $x^{k+1}=x^{k}+t_{k}d(x^{k}).$
	\item[Step 4] Compute $p(x^{k+1})$ and $\theta(x^{k+1})$, set $k\leftarrow k+1$, and go to \textbf{Step 1}.
\end{newlist}

In the step 3 of the CondG algorithm, we use the adaptative step size (see \cite{assunccao2021conditional}) to obtain $t_{k}$, that is,
$$t_{k}=\min\left\{1,\frac{\lvert\theta(x^k)\rvert}{L\|p(x^k)-x^k\|^2}\right\}.$$
Since $\theta(x)<0$ and $p(x)\neq x$ for non-Pareto critical points, the adaptative step size for 
the CondG algorithm is well-defined. The algorithm successfully stops if a Pareto critical point is found. Thus, hereafter, we assume that $\theta(x^{k})<0$ for all $k\geq0$, which means that the algorithm generates an infinite sequence $\{x^{k}\}$.

\section{Convergence analysis}\label{sec:4}

The following lemma indicates  that $\{x^{k}\}$ satisfies an important inequality, which can be proven similarly  to \cite[Proposition 13]{assunccao2021conditional}. It is noteworthy that a similar result has been further refined in our previous work \cite[Lemma 3]{chen2023conditional}.

\begin{lemma}\label{descent}
	For all $k\geq0$, it holds that
	\begin{equation}\label{d00}
		F(x^{k+1})-F(x^{k})\preceq-\dfrac{1}{2}\min\left\{\frac{\theta(x^{k})^{2}}{L\|p(x^{k})-x^{k}\|^{2}},-\theta(x^{k})\right\}e.
	\end{equation}
\end{lemma}

\begin{theorem}\label{convergence}
	Every limit point $\bar{x}$ of $\{x^{k}\}$ is a Pareto critical point of \eqref{mop}.
\end{theorem}
\textbf{Proof.}
	Let $\bar{x}\in\Omega$ be a limit point of $\{x_{k}\}$ and $\{x^{k_{j}}\}$ be a subsequence of $\{x_{k}\}$ such that $\lim_{j\rightarrow\infty}x^{k_{j}}=\bar{x}$. By the continuity argument of $F$, we have $\lim_{j\rightarrow\infty}F(x^{k_{j}})=F(\bar{x})$. Since $\{F(x^{k})\}$ is monotone decreasing as in Lemma \ref{descent}, it follows that $\lim_{k\rightarrow\infty}F(x^{k})=F(\bar{x})$, and thus
	\begin{equation}\label{convergence1}
		\lim_{k\rightarrow\infty}(F(x^{k+1})-F(x^{k}))=0.
	\end{equation}
	From the boundedness of $\{x^{k_{j}}\}$, and observing that Proposition \ref{prop2}, we know that $\{p(x^{k_{j}})\}$ is bounded. Let $\{p(x^{k_{j_{l}}})\}$ be a subsequence of $\{p(x^{k_{j}})\}$ such that $\lim_{l\rightarrow\infty}p(x^{k_{j_l}})=\bar{p}$. Consider the following two cases:
	
	\emph{Case 1.} Let $\bar{p}=\bar{x}$. By the definition of $\theta$ in \eqref{theta_xk} and the continuity argument of $JF$, we have
	\begin{equation*}
		\begin{aligned}
			\lim_{l\rightarrow\infty}\max_{i\in\langle m\rangle}\{\langle \nabla F_{i}(x^{k_{j_{l}}}),p^{k_{j_l}}-x^{k_{j_l}})\rangle\}
			=\max_{i\in\langle m\rangle}\lim_{l\rightarrow\infty}\{\langle \nabla F_{i}(x^{k_{j_{l}}}),p^{k_{j_l}}-x^{k_{j_l}})\rangle\}
			=\max_{i\in\langle m\rangle}\{\langle \nabla F_{i}(\bar{x}),\bar{p}-\bar{x}\rangle\}=0
		\end{aligned}
	\end{equation*}
	
	\emph{Case 2.} Let $\bar{p}\neq\bar{x}$. Combining \eqref{d00} with \eqref{convergence1}, we get
	\begin{equation*}
		\lim_{l\rightarrow\infty}	\min\left\{\frac{\theta(x^{k_{j_l}})^{2}}{L\|p(x^{k_{j_l}})-x^{k_{j_l}}\|^{2}},\lvert\theta(x^{k_{j_l}})\rvert\right\}=0.
	\end{equation*}
	It is clear that $\lim_{l\rightarrow\infty}\|p(x^{k_{j_l}})-x^{k_{j_l}}\|=\|\bar{p}-\bar{x}\|\neq0$. Therefore, $\lim_{l\rightarrow\infty}\theta(x^{k_{j_l}})=0$. According to \eqref{theta_xk}, we have
	\begin{equation}\label{convergence2}
		\theta(x^{k_{j_l}})\leq\max_{i\in\langle m\rangle}\{\langle\nabla F_{i}(x^{k_{j_l}}),u-x^{k_{j_l}}\rangle\}
	\end{equation}
	for all $u\in\Omega$. Taking the limit as $l\rightarrow\infty$ in \eqref{convergence2}, we have $\max_{i\in\langle m\rangle}\{\langle\nabla F_{i}(\bar{x}),u-\bar{x}\rangle\}\geq0$, which coincides with \eqref{opt_con}, and thus $\bar{x}$ is a Pareto critical point of \eqref{mop}.
\qed

\begin{remark}\normalfont
	In the proof of Theorem \ref{convergence}, we did not utilize the continuity of the function $\theta$ in \eqref{theta_xk}, which differs from the work in \cite[Remark 2]{assunccao2021conditional}.
\end{remark}

It follows from Lemma \ref{global} and Theorem \ref{convergence} that the following result holds.
\begin{theorem}
	If $F$ is convex on $\Omega$, then $\{x^{k}\}$ converges to a weak Pareto solution of  \eqref{mop}.
\end{theorem}

According to the definition of Pareto optimal solution and the process of descent methods in multiobjective optimization, the limit $$\lim_{k\rightarrow\infty}\min_{i\in\langle m
	\rangle}\{F_{i}(x^{k})-F_{i}(\bar{x})\}$$ indicates the convergence of the objectives, as reported in \cite{zeng2019convergence}. Actually,
the least reduction of the function values equals to zero in a descent method means that all
objective functions cannot decrease anymore. Next we give a result on the convergence rate
of $\{\min_{i\in\langle m
	\rangle}\{F_{i}(x^{k})-F_{i}(\bar{x})\}\}$. For  simplicity, let us define the following two constants:
\begin{equation}\label{gamma}
		\rho=\max\left\{\max_{i\in\langle m\rangle}\|\nabla F_{i}(x^{k})\|:k\geq0\right\}
\quad{\rm and}\quad
		\beta=\min\left\{\dfrac{1}{2\rho\sigma},\dfrac{1}{2L\sigma^{2}}\right\},
\end{equation}
where $\sigma=\sup\{\|p(x^{k})-x^{k}\|,k\geq0\}$.

\begin{theorem}
	If $F$ is convex on $\Omega$,
then
	\begin{equation}\label{rate00}
		\min_{i\in\langle m\rangle}\{F_{i}(x^{k})-F_{i}(\bar{x})\}\leq\dfrac{1}{\beta k}.
	\end{equation}
\end{theorem}
\textbf{Proof.}
	From Lemma \ref{descent}, and observing that $\theta(x^{k})<0$, for all $i\in\langle m\rangle$, we have
	\begin{equation}\label{rate01}
		F_{i}(x^{k})-F_{i}(x^{k+1})\geq\theta(x^{k})^{2}\min\left\{\frac{1}{2L\sigma^{2}},\dfrac{1}{2\lvert\theta(x^{k})\rvert}\right\}.
	\end{equation}
	According to \eqref{theta_xk} and the Cauchy--Schwarz inequality, it holds that
	\begin{equation*}
		\begin{aligned}
			\lvert\theta(x^{k})\rvert=\left\lvert\max_{i\in\langle m\rangle}\{\langle\nabla F_{i}(x^{k}),p(x^{k})-x^{k}\rangle\}\right\lvert
			\leq\max_{i\in\langle m\rangle}\{\|\nabla F_{i}(x^{k})\|\}\|p(x^{k})-x^{k}\|
			\leq\rho\sigma,
		\end{aligned}
	\end{equation*}
	which together with \eqref{gamma} and \eqref{rate01} gives us
	$
		F_{i}(x^{k})-F_{i}(x^{k+1})\geq\beta\theta(x^{k})^{2}
	$
	for all $i\in\langle m\rangle$. Therefore,
	\begin{equation*}
		F_{i}(x^{k})-F_{i}(\bar{x})\geq F_{i}(x^{k+1})-F_{i}(\bar{x})+\beta\theta(x^{k})^{2},
	\end{equation*}
	for all $i\in\langle m\rangle$. Taking the min with respect to $i\in\langle m\rangle$ on both sides of the above inequality, we have
	\begin{equation}\label{rate02}
		\min_{i\in\langle m\rangle}\{F_{i}(x^{k})-F_{i}(\bar{x})\}-\min_{i\in\langle m\rangle}\{F_{i}(x^{k+1})-F_{i}(\bar{x})\}\geq\beta\theta(x^{k})^{2}.
	\end{equation}
	Since $F$ is convex on $\Omega$, we get
	$
	F_{i}(\bar{x})- F_{i}(x^{k})\geq\langle \nabla F_{i}(x^{k}),\bar{x}-x^{k}\rangle
	$
	for all $i\in\langle m\rangle$, which combined with the relation \eqref{theta_xk} yields
	\begin{equation}\label{rate03}
		\begin{aligned}
			\max_{i\in\langle m\rangle}\{F_{i}(\bar{x})- F_{i}(x^{k})\}\geq\max_{i\in\langle m\rangle}\{\langle \nabla F_{i}(x^{k}),\bar{x}-x^{k}\rangle\}
			\geq\max_{i\in\langle m\rangle}\{\langle \nabla F_{i}(x^{k}),p(x^{k})-x^{k}\rangle\}
			=\theta(x^{k}).
		\end{aligned}
	\end{equation}
	According to Lemma \ref{descent}, we have  $F_{i}(\bar{x})\leq F_{i}(x^{k})$ for all $i\in\langle m\rangle$. Combing this with \eqref{rate03}, we get 
	$
			0\leq\min_{i\in\langle m\rangle}\{F_{i}(x^{k})-F_{i}(\bar{x})\}\leq-\theta(x^{k}),
	$
	and thus 
	\begin{equation}\label{rate04}
		\left(\min_{i\in\langle m\rangle}\{F_{i}(x^{k})-F_{i}(\bar{x})\}\right)^{2}\leq\theta(x^{k})^{2}.
	\end{equation}
	Let $a_{k}=\min_{i\in\langle m\rangle}\{F_{i}(x^{k})-F_{i}(\bar{x})\}$. Then, by \eqref{rate02} and \eqref{rate04}, we have
	$
		a_{k}-a_{k+1}\geq\beta a_{k}^{2}.
	$
	Thus, \eqref{rate00} follows immediately from Lemma \ref{sub_rate}.
\qed

\section{Numerical examples}\label{sec:5}

In this section, we present the numerical results of our  method to solve two multiobjective optimization problems with the unbounded feasible region.

\begin{example}\normalfont\label{ex1}
	Consider \eqref{mop} with $n=2$, $m=2$, 
	$
			F_{1}(x)=x_{1}+0.01(x_{2}+0.5)^{2}, F_{2}(x)=0.01(x_{1}+0.5)^{2}+x_{2}
	$
	and 
	$\Omega=\{x=(x_{1},x_{2})\in\mathbb{R}_{+}^{2}:x_{1}+x_{2}\geq1,x_{2}\geq0.5\}.$
	Both functions are convex on $\Omega$. Clearly,  $(\Omega^{\infty})^{*}=\Omega^{\infty}=\mathbb{R}_{+}^{2}$ and (A1) holds. 
\end{example} 

\begin{example}\normalfont\label{ex2}
	Consider \eqref{mop} with  $n=2$, $m=2$,
	$
			F_{1}(x)=-x_{1}+2x_{2}, F_{2}(x)=x_{1}+0.5\sin(x_{2})+1.1x_{2}
	$
	and 
	$\Omega=\{x=(x_{1},x_{2})\in\mathbb{R}^{2}:0.5x_{1}-x_{2}\leq0,-0.5x_{1}-x_{2}\leq0\}. $
	$F_{1}$ is convex on $\Omega$, whereas $F_{2}$ is not. Clearly, $(\Omega^{\infty})^{*}=\Omega^{\infty}=\Omega$ and (A1) holds. 
\end{example} 

\begin{figure}[htbp]
	\centering   
	\subfigure[Example \ref{ex1}] 
	{
		\begin{minipage}[b]{.45\linewidth} 
			\centering
			\includegraphics[scale=0.5]{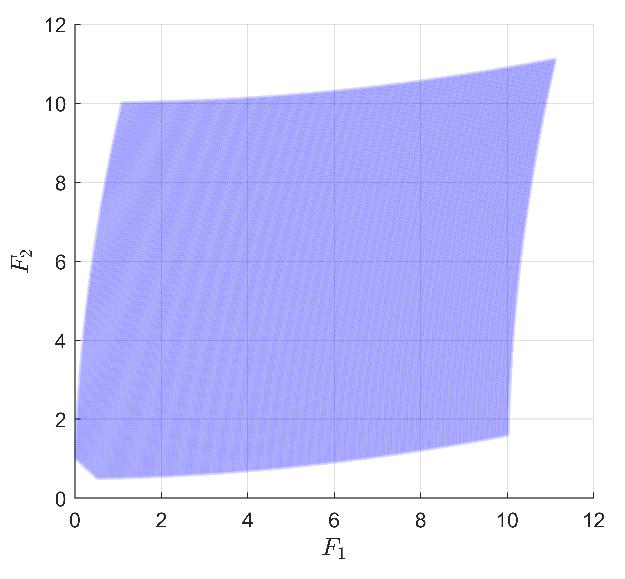}
		\end{minipage}
	}
	\subfigure[Example \ref{ex2}]
	{
		\begin{minipage}[b]{.45\linewidth}
			\centering
			\includegraphics[scale=0.5]{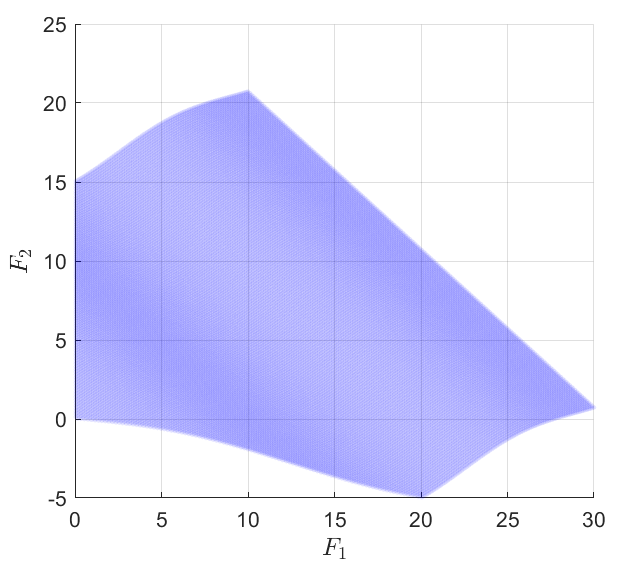}
		\end{minipage}
	}
	\caption{Visualization of the objective functions $F$ on Examples \ref{ex1} and \ref{ex2}.}
	\label{fig:pro1f}
\end{figure}

According to \cite[pp. 745]{assunccao2021conditional}, \eqref{sca_pro} is equivalent to the following optimization problem:
\begin{equation}\label{subpro}
	\begin{aligned}
		\text{min}\quad &\gamma\\
		\text{s.t.}\quad&\langle\nabla F_{i}(x),u-x\rangle\leq\gamma,~i\in\langle m\rangle,\\
		&u\in \Omega.
	\end{aligned}
\end{equation}
The experiments were conducted using MATLAB R2020b software on a PC with the following specifications: Intel i7-10700 processor running at 2.90 GHz and 32.00 GB RAM. The solver \texttt{fmincon} was employed to solve the subproblem \eqref{subpro}. The termination criterion (Step 1 of the CondG algorithm) was set as $\lvert\theta(x^{k})\rvert\leq\epsilon$ with $\epsilon=10^{-6}$. The maximum allowed number of outer iterations was set to 1000. For each test problem, the algorithm was run 100 times with initial points generated from a uniform random distribution within the respective feasible region.

Table \ref{res} presents the results obtained by the algorithm, organized into columns labeled ``it", ``gE", ``T" and ``\%." The ``it" column represents the average number of iterations, while ``gE" stands for the average number of gradient evaluations. The ``T" column indicates the average computational time (in seconds) to reach the critical point from an initial point, and ``\%" indicates the percentage of runs that have reached a critical point. As observed in Table \ref{res}, the algorithm can effectively solve the two given problems.

\begin{table}[htbp]
	\centering
	\caption{Performance of the algorithm on the two problems.}
	\vskip-0.1in
	\begin{tabular}{lllll}
		\hline
		& it    & gE    & T     & \% \\
		\hline
		Example 1 & 18.78  & 19.78  & 0.05 & 100 \\
		Example 2 & 14.81 & 15.81 & 0.04  & 100 \\
		\hline
	\end{tabular}%
	\label{res}%
\end{table}%

\begin{figure}[htbp]
	\centering   
	\subfigure[Example \ref{ex1}] 
	{
		\begin{minipage}[b]{.45\linewidth} 
			\centering
			\includegraphics[scale=0.5]{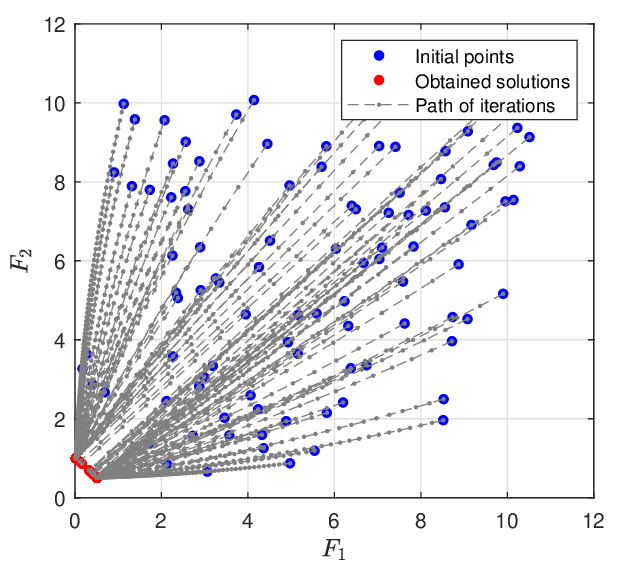}
		\end{minipage}
	}
	\subfigure[Example \ref{ex2}]
	{
		\begin{minipage}[b]{.45\linewidth}
			\centering
			\includegraphics[scale=0.5]{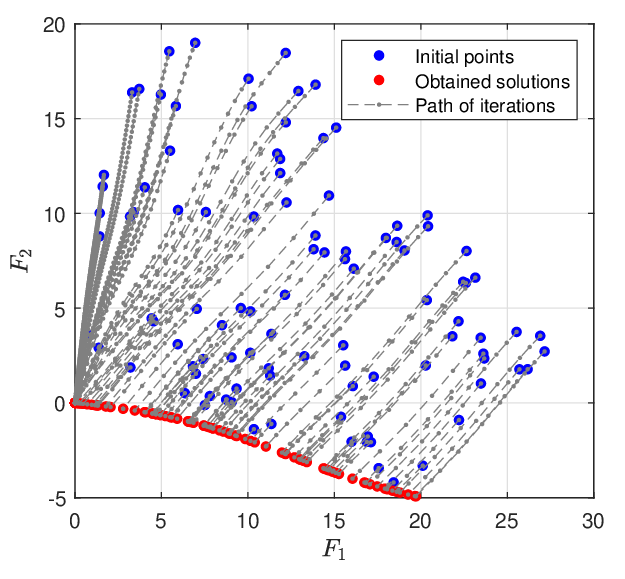}
		\end{minipage}
	}
	\caption{The final solutions and the paths of iterations obtained by the algorithm on the two examples.}
	\label{trajectories}
\end{figure}

To observe the movement of iteration points, we depict the trajectories of these points in Fig. \ref{trajectories}. In this figure, dashed lines represent the paths of algorithm iterations, blue points are the initial points, and red points correspond to the solutions found by the algorithm.

\end{document}